\documentclass[11pt,a4paper]{amsart}

\usepackage[utf8x]{inputenc} 
\usepackage[a4paper]{geometry} 
\usepackage{enumerate} 
\usepackage{amsmath} 
\usepackage{amssymb}
\usepackage{color,colortbl}
\usepackage[dvipsnames,table]{xcolor}
\usepackage[most]{tcolorbox}
\usepackage{fancyvrb}   
 
\definecolor{airforceblue}{rgb}{0.36, 0.54, 0.66}
\definecolor{bleudefrance}{rgb}{0.19, 0.55, 0.91}
\definecolor{darkorchid}{rgb}{0.6, 0.2, 0.8}
\definecolor{darkorange}{rgb}{1.0, 0.55, 0.0}
\definecolor{darkspringgreen}{rgb}{0.09, 0.45, 0.27}
\definecolor{commentoutput}{rgb}{0.40, 0.00, 0.0}
\definecolor{output}{rgb}{0.8, 0.0, 0.0}  
\definecolor{circOut}{rgb}{0.4, 1.0, 0.0} 

\definecolor{Gray}{gray}{0.7}  

\usepackage[arrow,curve,matrix,tips,frame,color]{xy}
\usepackage{adjustbox}
\usepackage{multirow}

\usepackage{verbatim}
\usepackage{hyperref}
                 \hypersetup{ pdfborder={0 0 0}, 
                              colorlinks=true, 
                              linkcolor=red, 
                              urlcolor=red,
                              citecolor=red,
                              linktoc=page,
                              pdfauthor={Giovanni Staglian{\`o}}, 
                              pdftitle={A package for computations with sparse resultants}
                            }                         

\theoremstyle{plain} 
\newtheorem{proposition}{Proposition}[section] 
\newtheorem{prop-def}[proposition]{Proposition-Definition}
\newtheorem{theorem}[proposition]{Theorem}

\theoremstyle{definition}

\newtheorem{example}[proposition]{Example} 
\theoremstyle{remark} 
\newtheorem{remark}[proposition]{Remark}

\newcommand{\PP}{{\mathbb{P}}}

\providecommand{\det}{\mathop{\rm det}\nolimits} 

\numberwithin{equation}{section}

\title{A package for computations with sparse resultants} 
\address{Dipartimento di Matematica e Informatica, Universit\`a degli Studi di Catania} 
\author[G. Staglian\`o]{Giovanni Staglian\`o}
\thanks{\emph{SparseResultants}, version 1.2, available at \url{https://faculty.math.illinois.edu/Macaulay2/doc/Macaulay2/share/doc/Macaulay2/SparseResultants/html/index.html}.} 
\email{\href{mailto:giovannistagliano@gmail.com}{giovannistagliano@gmail.com}} 
\keywords{Sparse resultant, sparse discriminant, hyperdeterminant} 
\subjclass[2020]{
   13P15, 
   68W30
}

\begin{document}

\begin{abstract} 
We introduce the \emph{Macaulay2} package 
\href{https://faculty.math.illinois.edu/Macaulay2/doc/Macaulay2/share/doc/Macaulay2/SparseResultants/html/index.html}{\emph{SparseResultants}},
which provides 
 general tools for computing sparse resultants, sparse discriminants, and 
hyperdeterminants. 
We give some background on the theory
and briefly 
show 
 how the package works.
\end{abstract}

\maketitle

\section*{Introduction}
The classical Macaulay 
resultant \cite{macaulaypaper1902} (also called the dense resultant) of a system of $n+1$ polynomial equations 
in $n$ variables 
characterizes the solvability of the system, and therefore
it
is a fundamental tool in computer algebra. 
However, it is a large polynomial, 
since it depends on all  coefficients of the equations.
If we restrict attention 
to \emph{sparse} polynomial equations, 
that is, to polynomials which involve only  monomials lying in a small set,
then we can replace the \emph{dense} resultant with the  \emph{sparse} resultant.

The sparse resultant generalizes not only the dense 
resultant but,
for specific choices of the set of monomials, 
we can obtain other types of classical resultants, such as for instance 
the \emph{Dixon resultant} \cite{Dixon1908} 
and the \emph{hyperdeterminant} \cite{cayley_2009,GELFAND1992226}. 
In the last decades, sparse resultants 
have received a lot of interest, both  
from a theoretical point of view (see \emph{e.g.}, \cite{gelfandbook,NewtonSturm,CattaniDickenstein1998,PoissonSombraAndrea})
and from more computational and applied aspects 
(see \emph{e.g.},  \cite{EMIRIS19993,CaEmi2000,Sturmfels02solvingsystems,DAndrea2001MacaulaySF,JeronimoSabia2004,coxbook,Matera2009Jeronimo}).

Using the computer program \emph{Macaulay2} \cite{macaulay2},
dense resultants can be calculated using the package \emph{Resultants} \cite{packageResultants},
while sparse resultants can be calculated using the new package \emph{SparseResultants}.
We point out that in the latter most of the algorithms implemented 
are based on elimination via Gr\"obner basis methods.
The main defect of this approach is that even when the input polynomials have numerical coefficients, in the calculation all the coefficients are replaced by variables. 
However, 
this approach suffices  for a number of applications, as we try to show in the following.

This short paper is organized as follows. In Section~\ref{sec 1}, we review the general theory
of sparse resultants (subsections~\ref{sparse-mixed} and \ref{sparse unmixed})  and related topics as the {sparse discriminants} (subsection~\ref{discriminants}) and the {hyperdeterminants} (subsection~\ref{sec dets}).
We focus  on the computational aspects used in the package \emph{SparseResultants}.
In Section~\ref{computations},
we illustrate how this package works 
with the help of some examples.

\section{An overview of sparse elimination}\label{sec 1}
In this section we give some background  on the theory of sparse resultants, 
sparse discriminants, and hyperdeterminants. For details and proofs we refer mainly to 
\cite[Chapters~8, 9, 13, and 14]{gelfandbook} and \cite[Chapter~7]{coxbook}; other references 
are \cite{SturmfelsElimination,OttSurvey}.
\subsection{Sparse mixed resultant}\label{sparse-mixed}
Let $R = \mathbb{C}[x_1^{\pm1},\ldots,x_n^{\pm1}]$ be 
the ring of complex Laurent polynomials in $n$ variables.
The set of monomials in $R$ is identified with $\mathbb{Z}^n$ 
by associating to $x^{\omega} = x_1^{\omega_1}\cdots x_n^{\omega_n}\in R$ 
the exponent vector $\omega = (\omega_1,\ldots,\omega_n)\in\mathbb{Z}^n$.
If $\mathcal A$ is a finite subset of $\mathbb{Z}^n$, we denote by 
$\mathbb{C}^{\mathcal A}$ 
the space of polynomials in $R$
involving only monomials from $\mathcal A$, that is, of polynomials 
of the form
$\sum_{\omega\in \mathcal A} a_{\omega} x^{\omega}$.

Let $\mathcal A_0,\ldots,\mathcal A_n$ be $n+1$ finite subsets of $\mathbb{Z}^n$ 
satisfying the following conditions:
\begin{enumerate}
 \item each $\mathcal A_i$ generates $\mathbb{R}^n$ as an affine space; and
 \item the union of the sets  $\mathcal A_i$ generates 
 $\mathbb{Z}^n$ as a $\mathbb{Z}$-module.
\end{enumerate}
Let $\mathcal{Z}_{\mathcal A_0,\ldots,\mathcal A_n}\subset \prod_{i=0}^n \mathbb{C}^{\mathcal A_{i}}$ be 
the Zariski closure in the product $\prod_{i=0}^n \mathbb{C}^{\mathcal A_{i}}$ of the set 
\begin{equation}\label{formula sparse res}
 \{(f_0,\ldots,f_n)\in \prod_{i=0}^n \mathbb{C}^{\mathcal A_{i}}: \exists x\in(\mathbb C^{\ast})^{n} \mbox{ such that }f_0(x) = \cdots = f_n(x) = 0\},
\end{equation}
where $\mathbb{C}^{\ast} = \mathbb{C}\setminus\{0\}$ and $f_i(x)=\sum_{\omega\in \mathcal A_i} a_{i,\omega} x^{\omega}$, for $i=0,\ldots,n$.
\begin{prop-def}[\cite{gelfandbook}, Chapter~8, \S 1]
 Under the above assumptions, the variety $\mathcal{Z}_{\mathcal A_0,\ldots,\mathcal A_n}$ 
 is an irreducible hypersurface in $\prod_{i=0}^n \mathbb{C}^{\mathcal A_i}$ that can be
 defined by an integral irreducible polynomial $\mathrm{Res}_{\mathcal A_0,\ldots,\mathcal A_n}\in\mathbb{Z}[(a_{i,\omega}),i=0,\ldots,n]$ 
 in the coefficients $a_{i,\omega}$ of $f_{i}$, for $i=0,\ldots, n$.
 Such a polynomial $\mathrm{Res}_{\mathcal A_0,\ldots,\mathcal A_n}$ is unique 
 up to sign 
 and is called the $(\mathcal A_0,\ldots,\mathcal A_n)$-resultant (also known as the \emph{sparse (mixed) resultant}).
\end{prop-def}
The polynomial
$\mathrm{Res}_{\mathcal A_0,\ldots,\mathcal A_n}$ is homogeneous with respect to each group 
of variables $(a_{i,\omega})$, for $i=0,\ldots,n$. Moreover, $\mathrm{Res}_{\mathcal A_0,\ldots,\mathcal A_n}(f_0,\ldots,f_n) = 0$
if the $(n+1)$-tuple $(f_0,\ldots,f_n)$ belongs to \eqref{formula sparse res}.
\begin{example}
 Let $d_0,\ldots,d_n$ be positive integers. For $i=0,\ldots,n$, let 
 $\mathcal A_i = \{\omega=(\omega_1,\ldots,\omega_n)\in \mathbb{Z}_{\geq 0}^n: \sum_{j=1}^n \omega_j \leq d_i\}$.
 Then the $(\mathcal A_0,\ldots,\mathcal A_n)$-resultant coincides with the classical (affine) resultant $\mathrm{Res}_{d_0,\ldots,d_n}$,
 also called the \emph{dense resultant}.
 Therefore, 
 if $F_i\in \mathbb{C}[x_0,x_1,\ldots,x_n]$ 
 denotes the polynomial
 obtained by homogenizing $f_i\in\mathbb{C}^{\mathcal A_i}$ with respect to a new variable $x_0$,
 then 
 $\mathrm{Res}_{\mathcal A_0,\ldots,\mathcal A_n}(f_0,\ldots,f_n) = 0 $
 if and only if $F_0,\ldots,F_n$ have a common non-trivial root.
\end{example}

\subsection{Sparse unmixed resultant}\label{sparse unmixed}
Keep the notation and assumptions as above. 
If all the sets $\mathcal A_i$ coincide with each other, that is, 
$\mathcal A_0 = \cdots = \mathcal A_n = \mathcal A$, then the $(\mathcal A_0,\ldots,\mathcal A_n)$-resultant 
is called the $\mathcal A$-resultant (also known as the \emph{sparse (unmixed) resultant}). 
In this  case, we have a useful geometric interpretation that
allows
us to write down the $\mathcal{A}$-resultant in a compact form.
Indeed, by choosing a numbering $\omega^{(0)},\ldots,\omega^{(k)}$ of the elements of $\mathcal A$, 
we get a map $\phi_{\mathcal A}:(\mathbb{C}^{\ast})^n\to \PP^{k}$ 
defined by $\phi_{\mathcal A}(x) = (\omega^{(0)}(x):\cdots:\omega^{(k)}(x))$.
Let $X_{\mathcal A}\subset \PP^k$ be the closure of the image of $\phi_{\mathcal A}$,
which is an irreducible toric variety of dimension $n$.
Then, by taking pull-backs
we get an identification 
between 
the space of polynomials in $\mathbb{C}^{\mathcal A}$ 
with the space of linear forms on $\PP^k$. Moreover,
if $f_0,\ldots,f_n\in \mathbb{C}^{\mathcal A}$ have a common root in $(\mathbb{C}^{\ast})^n$
then the corresponding linear forms $l_0,\ldots,l_n$ on $\PP^k$ define 
a linear subspace that intersects $X_{\mathcal A}$. From this, 
the following proposition follows directly.
\begin{proposition}[\cite{gelfandbook}, Chapter~8, \S 2]\label{chowUnmixed}
 The polynomial $\mathrm{Res}_{\mathcal A}\in\mathbb{Z}[a_0^{(i)},\ldots,a_k^{(i)},i = 0,\ldots ,n]$ 
 coincides with the $X$-resultant of $X_{\mathcal A}\subset\PP^k$. More precisely, 
 let $W_{\mathcal A}\subset \mathbb{G}(k-n-1,\PP^k)$ be the Chow hypersurface
 of the variety $X_{\mathcal A}$, 
 and let \[\psi:\PP(\mathbb{C}^{(n+1)\times (k+1)})\dashrightarrow \mathbb{G}(n,k)\simeq \mathbb{G}(k-n-1,k)\]
 be the natural projection from the projectivization of the space of  complex matrices of shape $(n+1)\times (k+1)$
 to $\mathbb{G}(n,k)$.
 Then we have that 
  $\mathrm{Res}_{\mathcal A}$ 
 is the polynomial defining the pull-back $\overline{\psi^{-1}(W_{\mathcal A})}$.
\end{proposition}
\begin{remark}
 With the notation of the proposition above,  in coordinates 
the map $\psi$
 is defined by the $(n+1)\times (n+1)$ minors 
 of the generic $(n+1)\times (k+1)$  matrix of variables 
 \begin{equation}\label{matrix}
  \begin{pmatrix}
     a_0^{(0)}& a_1^{(0)}&\cdots&a_k^{(0)}\\
    \vdots & \vdots & \ddots &\vdots \\
     a_0^{(n)}& a_1^{(n)}&\cdots&a_k^{(n)}\\
    \end{pmatrix} .
 \end{equation}
In particular,
 $\mathrm{Res}_{\mathcal A}$ can be expressed as a homogeneous polynomial of degree $\deg(X_{\mathcal{A}})$ in the $(n+1)\times (n+1)$ minors 
 of the  matrix \eqref{matrix}. 
\end{remark}
\begin{example}
 Let $\mathcal{A} = \{(\omega_1,\omega_2)\in\mathbb{Z}^2:\omega_1+\omega_2\leq 2\}$,
 so that $X_{\mathcal{A}}\subset\PP^{5}$ is the Veronese surface.
 The $\mathcal{A}$-resultant 
 is a polynomial of degree $12$ in $18$ variables with $21,894$ terms. It can be expressed as a polynomial 
 of degree $4$ in the Pl\"ucker coordinates of $\mathbb{G}(2,5)$ with $74$ terms.
\end{example}
\subsection{Sparse discriminant}\label{discriminants}
We continue by letting $\mathcal A\subset\mathbb{Z}^n$  be a finite set of 
 $k+1$ elements that generate $\mathbb{Z}^n$ as a $\mathbb{Z}$-module,
 and let $\phi_{\mathcal A}:(\mathbb C^{\ast})^n\to\mathbb{P}^k$ and $X_{\mathcal A}\subset\PP^{k}$ be 
 defined as above. 
 Let $\nabla_{\mathcal A}\subset\mathbb{C}^{\mathcal A}$ be the Zariski closure of the set 
 \begin{equation}\label{disc set}
  \{f\in\mathbb{C}^{\mathcal A} : \exists x\in (\mathbb{C}^{\ast})^n \mbox{ such that } f(x) = \frac{\partial f}{\partial x_1}(x)=\cdots=\frac{\partial f}{\partial x_n}(x)=0\}.
 \end{equation}
\begin{prop-def}[\cite{gelfandbook}, Chapter~9, \S 1]
 The projectivization $\mathbb{P}(\nabla_{\mathcal A})\subset\PP^{k}$ of the variety $\nabla_{\mathcal A}$ 
 coincides with the  dual variety $X_{\mathcal A}^{\vee}$ of $X_{\mathcal A}$.
 In the case where $X_{\mathcal A}^{\vee}$ is a hypersurface, 
 an integral irreducible polynomial $\mathrm{Disc}_{\mathcal A}$ defining it (which is unique up to sign) 
 is called the $\mathcal A$-discriminant (also known as the \emph{sparse discriminant}).
\end{prop-def}
Thus the $\mathcal{A}$-discriminant  (when it exists) is 
a
homogeneous polynomial 
$\mathrm{Disc}_{\mathcal A}\in \mathbb{Z}[a_{\omega}, \omega\in \mathcal A]$, 
and $\mathrm{Disc}_{\mathcal A}(f)=0$ for each polynomial $f$ belonging to \eqref{disc set}.
\begin{example}
 Let $d\geq 1$ and let 
 $\mathcal A = \{(\omega_1,\ldots,\omega_n)\in \mathbb{Z}_{\geq 0}^n: \sum_{j=1}^n \omega_j \leq d\}$.
 Then the $\mathcal A$-discriminant coincides with the classical (affine) discriminant $\mathrm{Disc}_{d}$,
 also called the \emph{dense discriminant}.
 Therefore, 
 if $F\in \mathbb{C}[x_0,x_1,\ldots,x_n]$ 
 denotes the polynomial
 obtained by homogenizing $f\in\mathbb{C}^{\mathcal A}$ with respect to a new variable $x_0$,
 then 
 $\mathrm{Disc}_{\mathcal A}(f) = 0 $
 if and only if the hypersurface $\{F = 0\}\subset\PP^{n}$ is not smooth.
\end{example}
\begin{remark}[``Cayley trick'', \cite{gelfandbook}, Chapter~9, Proposition~1.7]
 Let $\mathcal A_0,\ldots,\mathcal A_n\subset\mathbb{Z}^n$ be finite subsets 
 satisfying the assumptions of subsection~\ref{sparse-mixed}.
 Let $\mathcal A\subset\mathbb{Z}^n\times \mathbb{Z}^n$ be defined by
 $\mathcal A = (\mathcal A_0\times \{0\})\cup (\mathcal A_1\times \{e_1\}) \cup \cdots \cup (\mathcal A_n\times \{e_n\})$,
 where the $e_i$ are the standard basis vectors of $\mathbb{Z}^n$. Thus a polynomial 
 $f\in \mathbb{C}^{\mathcal A}$ has the form 
 \[
  f_0(x) + \sum_{i = 1}^n y_i\,f_i(x) \in \mathbb{C}[x_1,\ldots,x_n,y_1,\ldots,y_n],
 \]
where $f_i\in \mathbb{C}^{\mathcal A_i}$.
We have the following relation (up to sign), known as the ``Cayley trick'': 
\begin{equation}\label{cayleyForm}
  \mathrm{Res}_{\mathcal A_0,\ldots,\mathcal A_n}(f_0,\ldots,f_n) = \mathrm{Disc}_{\mathcal A}(f_0(x) + \sum_{i = 1}^n y_i\,f_i(x)).
\end{equation}
\end{remark}

\subsection{Hyperdeterminant}\label{sec dets}
An important special type of sparse discriminant is the 
determinant (or hyperdeterminant) of multidimensional matrices,
which was introduced by Cayley \cite{cayley_2009} (see also \cite[Chapter~14]{gelfandbook} and \cite{OttSurvey}).
Let $f$ be a multilinear form 
in $r$ groups of variables 
$x_0^{(1)},\ldots,x_{k_1}^{(1)};\ldots ; x_0^{(r)},\ldots,x_{k_r}^{(r)}$, that is 
\[f = \sum_{\begin{subarray}{c} 0\leq i_{\iota}\leq k_{\iota}  \end{subarray}} a_{i_1,\ldots,i_r} x_{i_1}^{(1)}\cdots x_{i_r}^{(r)}.\]
Let $\mathcal{A}\subset\mathbb{Z}^{(k_1+1)+\cdots+(k_r+1)}$ denote the set 
of exponent vectors that can occur in such a form $f$.
Notice that to give $f$ is equivalent to giving 
an
$r$-dimensional matrix
 \[M_f = (a_{i_1,\ldots,i_r})_{0\leq i_{\iota}\leq k_{\iota}}\]
 of shape $(k_1+1)\times \cdots \times (k_r+1)$.
 The determinant of shape $(k_1+1)\times\cdots\times (k_r+1)$ is defined to be the 
 $\mathcal{A}$-discriminant, that is, for a form $f$ as above, we have
 \[
 \det(M_f) = \mathrm{Disc}_{\mathcal{A}}(f) .
 \]
One sees that the variety $X_{\mathcal A}$ 
is the image of the Segre embedding of $\PP^{k_1}\times\cdots\times\PP^{k_r}$.
Therefore,  
the hypersurface 
in $\PP(\mathbb{C}^{(k_1+1)\times\cdots\times(k_r+1)})$ defined by 
the determinant of shape $(k_1+1)\times\cdots\times(k_r+1)$
is the dual variety of $\PP^{k_1}\times\cdots\times\PP^{k_r}$.
Notice also that
we have $\det(M_f) = 0$ if and only if the hypersurface $\{f = 0\}\subset \PP^{k_1}\times\cdots\times\PP^{k_r}$ is not smooth.

The next two basic results have been proved in  \cite[Chapter~14, Theorems~1.3 and 2.4]{gelfandbook}.
\begin{theorem}[\cite{gelfandbook}]
 The determinant of shape $(k_1+1)\times\cdots\times(k_r+1)$ 
 exists (that is the dual variety of $\PP^{k_1}\times\cdots\times\PP^{k_r}$ is a hypersurface) if and only if 
 \begin{equation}\label{cond}
  2\,\max_{{1\leq j\leq r}}(k_j) \leq \sum_{j=1}^r k_j .
 \end{equation}
 \end{theorem}
\begin{theorem}[\cite{gelfandbook}]
 Denote by $N(k_1,\ldots,k_r)$ the 
 degree of the determinant of shape $(k_1+1)\times\cdots\times(k_r+1)$
 when \eqref{cond} is satisfied,  and let $N(k_1,\ldots,k_r) = 0$ otherwise.
 We have 
 $$
 \sum_{k_1,\ldots,k_r\geq 0} N(k_1,\ldots,k_r)z_1^{k_1}\cdots z_r^{k_r} 
 = \frac{1}{(1-\sum_{i=2}^r(i-2)e_i(z_1,\ldots,z_r))^2},
 $$
 where $e_i(z_1,\ldots,z_r)$ is the $i$-th elementary symmetric polynomial. 
\end{theorem}
\begin{remark}[\cite{gelfandbook}, Chapter~4, Propositions~1.4 and 1.8]
The determinant of shape $(k_1+1)\times\cdots\times (k_r+1)$
is invariant under 
the action of $\mathrm{SL}(k_1+1)\times\cdots\times\mathrm{SL}(k_r+1)$
on the space of matrices of shape $(k_1+1)\times\cdots\times(k_r+1)$. It is also invariant 
under permutations of the dimensions, that is,
if $M = (a_{i_1,\ldots,i_r})$ is a matrix of shape $(k_1+1)\times\cdots\times (k_r+1)$
and $\sigma$ is a permutation of indices $1,\ldots,r$,
denoting by $\sigma(M)$ the matrix of shape $(k_{\sigma^{-1}(1)}+1)\times\cdots\times(k_{\sigma^{-1}(r)}+1)$,
whose $(i_1,\ldots,i_r)$-th entry is equal to $a_{i_{\sigma(1)},\ldots,i_{\sigma(r)}}$,
we have $\det(\sigma(M)) = \det(M)$.
\end{remark}
There are at least two important cases where
determinants can be computed 
without resorting to elimination. We briefly recall them in \ref{schlafli-sec} and \ref{boundary}.

\subsubsection{Schl\"afli's method}\label{schlafli-sec}
Let $M$ be an $r$-dimensional matrix of shape 
$(k_1+1)\times\cdots\times (k_r+1)$ 
corresponding to a multilinear form 
$f\in \mathbb{C}[x_0^{(1)},\ldots,x_{k_1}^{(1)};\ldots ; x_0^{(r)},\ldots,x_{k_r}^{(r)}]$. Assume that 
there exist both the determinants 
of shapes $(k_1+1)\times\cdots\times(k_r+1)$ and
$(k_1+1)\times\cdots\times(k_{r-1}+1)$.
We can interpret the $r$-dimensional matrix $M$ as an 
$(r-1)$-dimensional matrix $\tilde{M}(x_0^{(r)},\ldots,x_{k_r}^{(r)})$
of shape $(k_1+1)\times \cdots\times (k_{r-1}+1)$ whose entries 
are linear forms in the variables $x_0^{(r)},\ldots,x_{k_r}^{(r)}$;
in other words, we can see $f$ 
as a polynomial 
$\tilde{f}\in (\mathbb{C}[x_0^{(r)},\ldots,x_{k_r}^{(r)}])[x_0^{(1)},\ldots,x_{k_1}^{(1)};\ldots ; x_0^{(r-1)},\ldots,x_{k_{r-1}}^{(r-1)}]$.
Let $F_M = F_M(x_0^{(r)},\ldots,x_{k_r}^{(r)}) = \det(\tilde{M}(x_0^{(r)},\ldots,x_{k_r}^{(r)}))$,
which is a homogeneous polynomial in $x_0^{(r)},\ldots,x_{k_r}^{(r)}$,
and let $\mathrm{Disc}(F_M)$ be the (classical) discriminant of $F_M$.
Then we have the following:
\begin{theorem}[\cite{gelfandbook,bhlpart219634}]
 The polynomial $\mathrm{Disc}(F_M)$ is divisible by the determinant $\det(M)$.
 Moreover if the shape of $M$ is one of the following 
 \begin{equation}\label{sch}
 m\times m\times 2,\ m\times m\times 3,\ 2\times2\times 2 \times 2, \mbox{ with } m\geq2,
    \end{equation}
then we have $\mathrm{Disc}(F_M) = \det(M)$.
    \end{theorem}
The  method above turns out to be very effective; however
it was conjectured in \cite[p.~479]{gelfandbook}, and later proved in \cite{AIF19964635910},
that the shapes in \eqref{sch} are the only ones 
for which the method gives the determinant exactly.

\subsubsection{Determinants of boundary shape}\label{boundary}
For an $(r+1)$-dimensional matrix $M$ of shape $(k_0+1)\times (k_1+1)\times \cdots\times(k_r+1)$,
we say that it is of \emph{boundary shape} if the inequality 
\eqref{cond} is an equality. \emph{W.l.o.g.}, we can assume that $k_0 = \max_{0\leq j\leq r} (k_j)$,
so that $k_0 = k_1+\cdots+k_{r}$.
Let $f\in\mathbb{C}[x_0^{(0)},\ldots,x_{k_0}^{(0)};\ldots ; x_0^{(r)},\ldots,x_{k_r}^{(r)}]$ 
be the corresponding multilinear form of such a matrix $M$.
Thinking $f$ as a linear polynomial in 
$(\mathbb{C}[x_0^{(1)},\ldots,x_{k_1}^{(1)};\ldots ; x_0^{(r)},\ldots,x_{k_{r}}^{(r)}])[x_0^{(0)},\ldots,x_{k_0}^{(0)}]$, we can interpret $M$ 
as a list of $k_0+1 $ 
multilinear forms $f_0,\ldots,f_{k_0}$ in the $r$ groups of variables
$x_0^{(1)},\ldots,x_{k_1}^{(1)};\ldots ; x_0^{(r)},\ldots,x_{k_{r}}^{(r)}$.

A simple consequence of the ``Cayley trick'' (see \cite[Chapter~3, Corollary~2.8]{gelfandbook}) gives the following:
\begin{proposition}[\cite{gelfandbook}]\label{boundary det}
 The determinant of an $(r+1)$-dimensional matrix $M$ of boundary shape $(k_0+1)\times\cdots\times(k_r+1)$ 
 coincides with the resultant of the multilinear forms $f_0,\ldots,f_{k_0}$,
 that is,
 $\det(M) = 0 $ if and only if the system of multilinear equations
 $f_0(x) = \cdots = f_{k_0}(x) = 0$ has a non-trivial solution on $\PP^{k_1}\times \cdots \times \PP^{k_{r}}$.
 In other words, the determinant of shape $(k_0+1)\times \cdots\times(k_r+1)$ 
 coincides 
 with the $X$-resultant of the Segre embedding of $\PP^{k_1}\times \cdots \times \PP^{k_{r}}$.
\end{proposition}
\begin{remark}\label{remBoundaryDet}
We also mention that the determinant of a matrix $M$ of boundary 
shape $(k_0+1)\times\cdots\times(k_r+1)$ 
can be explicitly expressed 
as the determinant of an ordinary square matrix of order $\frac{(k_0+1)!}{k_1!\cdots k_r!}$
whose entries are linear forms in the entries of $M$; see \cite[Chapter~14, Theorem~3.3]{gelfandbook}.
\end{remark}

\section{Sparse resultants in \emph{Macaulay2}}\label{computations}
In this section, we briefly describe 
some of the 
functions 
implemented in the package 
\href{https://faculty.math.illinois.edu/Macaulay2/doc/Macaulay2/share/doc/Macaulay2/SparseResultants/html/index.html}{\emph{SparseResultants}}. 
For more details and examples, we refer to its documentation. 

One of the main functions  is  \texttt{sparseResultant},
which via elimination techniques calculates 
sparse mixed resultants $\mathrm{Res}_{\mathcal{A}_0,\ldots,\mathcal{A}_n}$ (see subsection~\ref{sparse-mixed})
and sparse unmixed resultants $\mathrm{Res}_{\mathcal{A}}$ (see subsection~\ref{sparse unmixed}).
This function can be called in two ways.
The first one is to pass a list of $n+1$ 
matrices $A_0,\ldots,A_n$ over $\mathbb{Z}$ and with $n$ rows to represent the
sets $\mathcal{A}_0,\ldots,\mathcal{A}_n\subset \mathbb{Z}^n$ 
(it is enough to pass just one matrix $A$ in the unmixed case).
Then the output will be another 
 function that takes
 as input $n+1$ polynomials $f_i=\sum_{\omega\in\mathcal{A}_i} a_{i,\omega} x^{\omega}$, 
for $i = 0,\ldots,n$, and returns their sparse resultant. 
An error is thrown 
if the polynomials $f_i$ do not have the correct form.
Roughly, this returned function  
is a container for the general expression of the sparse resultant 
(possibly written out in a compact form as in Proposition~\ref{chowUnmixed})
and for the rule to evaluate it at the $n+1$ polynomials $f_i$.
The second way to call \texttt{sparseResultant} is to pass directly the polynomials $f_i$.
This is equivalent to forming the matrices  $A_i$ whose 
 columns  are given by
$\{\omega\in\mathbb{Z}^n:\mbox{ the coefficient in }f_i \mbox{ of }x^{\omega} \mbox{ is }\neq 0\}$ (see the function \texttt{exponentsMatrix}) 
and then proceeding as described above.

As an example we now calculate 
a particular type of sparse unmixed resultant, known as 
 the
\emph{Dixon resultant} (see \cite[Section~2.4]{SturmfelsElimination} and \cite[Chapter~7, \S 2, Exercise~10]{coxbook};
see also the classical reference \cite{Dixon1908}).
\begin{example}
Consider the following system of three bihomogeneous polynomials of bidegree $(2,1)$ in the two groups of variables $(x_0,x_1),(y_0,y_1)$:
\begin{align}
 c_{1,1} x_1^2 y_1+c_{1,2} x_1 x_2 y_1+c_{1,3} x_2^2 y_1+c_{1,4} x_1^2 y_2+c_{1,5} x_1 x_2 y_2+c_{1,6} x_2^2 y_2 &= 0,\nonumber \\
 c_{2,1} x_1^2 y_1+c_{2,2} x_1 x_2 y_1+c_{2,3} x_2^2 y_1+c_{2,4} x_1^2 y_2+c_{2,5} x_1 x_2 y_2+c_{2,6} x_2^2 y_2 &=0, \label{syst} \\
 c_{3,1} x_1^2 y_1+c_{3,2} x_1 x_2 y_1+c_{3,3} x_2^2 y_1+c_{3,4} x_1^2 y_2+c_{3,5} x_1 x_2 y_2+c_{3,6} x_2^2 y_2 &= 0.\nonumber
\end{align}
Putting $x_2 = y_2 = 1$ we get a system of three non-homogeneous polynomials in two variables $(x,y)=(x_1,y_1)$,
of which we can calculate the sparse (unmixed) resultant. This polynomial 
is homogeneous of degree $12$ in the $18$ variables $c_{1,1},\ldots,c_{3,6}$ with $20,791$ terms, which 
vanishes if and only if \eqref{syst} has a nontrivial solution. 
 The time for this computation is less than one second (on a standard laptop computer).
\begin{tcolorbox}[breakable=true,boxrule=0pt,opacityback=0.0,enhanced jigsaw]
{\footnotesize
\begin{Verbatim}[commandchars=&\{\}]
&colore{darkorange}{$ M2 --no-preload}
&colore{output}{Macaulay2, version 1.18}
&colore{darkorange}{i1 :} &colore{airforceblue}{needsPackage} "&colore{bleudefrance}{SparseResultants}";
&colore{darkorange}{i2 :} R = ZZ[c_(1,1)..c_(3,6)][x,y];
&colore{darkorange}{i3 :} f = (c_(1,1)*x^2*y+c_(1,2)*x*y+c_(1,3)*y+c_(1,4)*x^2+c_(1,5)*x+c_(1,6),
          c_(2,1)*x^2*y+c_(2,2)*x*y+c_(2,3)*y+c_(2,4)*x^2+c_(2,5)*x+c_(2,6),
          c_(3,1)*x^2*y+c_(3,2)*x*y+c_(3,3)*y+c_(3,4)*x^2+c_(3,5)*x+c_(3,6));
&colore{darkorange}{i4 :} A = &colore{bleudefrance}{exponentsMatrix} f
&colore{circOut}{o4 =} &colore{output}{| 0 0 1 1 2 2 |}
&colore{output}{     | 0 1 0 1 0 1 |}
&colore{output}{              2        6}
&colore{circOut}{o4 :} &colore{output}{Matrix ZZ  <--- ZZ}
&colore{darkorange}{i5 :} &colore{darkorchid}{time} Res = &colore{bleudefrance}{sparseResultant} A;
&colore{commentoutput}{     -- used 0.241391 seconds}
&colore{circOut}{o5 :} &colore{output}{SparseResultant (sparse unmixed resultant associated to | 0 0 1 1 2 2 |)}
&colore{output}{                                                             | 0 1 0 1 0 1 |}
&colore{darkorange}{i6 :} &colore{darkorchid}{time} U = Res f;
&colore{commentoutput}{     -- used 0.574002 seconds}
&colore{darkorange}{i7 :} (&colore{airforceblue}{first degree} U, # &colore{airforceblue}{terms} U)
&colore{circOut}{o7 =} &colore{output}{({12}, 20791)}
\end{Verbatim}
} 
\end{tcolorbox}
\end{example}

Another function of the package is \texttt{sparseDiscriminant},
which calculates 
sparse discriminants $\mathrm{Disc}_{\mathcal{A}}$ (see subsection~\ref{discriminants}).
This function works similar to the previous one. In particular, it accepts as input either 
a matrix representing the exponent vectors of a (Laurent) polynomial
or just directly the polynomial.
\begin{example}
Using the Cayley trick (formula \eqref{cayleyForm}), we
express the dense resultant of three generic ternary forms of degrees $1,1,2$ (which is a special type of sparse mixed resultant)
as a sparse discriminant. The calculation time is less than one second.
 \begin{tcolorbox}[breakable=true,boxrule=0pt,opacityback=0.0,enhanced jigsaw]
{\footnotesize
\begin{Verbatim}[commandchars=&\{\}]
&colore{darkorange}{i8 :} &colore{airforceblue}{clearAll};
&colore{darkorange}{i9 :} K = ZZ[a_0..a_2,b_0..b_2,c_0..c_5], Rx = K[x_1,x_2];
&colore{darkorange}{i10 :} f = (a_0+a_1*x_1+a_2*x_2, 
           b_0+b_1*x_1+b_2*x_2, 
           c_0+c_1*x_1+c_2*x_2+c_3*x_1^2+c_4*x_1*x_2+c_5*x_2^2);
&colore{darkorange}{i11 :} Rxy = K[x_1,x_2,y_1,y_2], f' = (&colore{airforceblue}{sub}(f_0,Rxy), &colore{airforceblue}{sub}(f_1,Rxy), &colore{airforceblue}{sub}(f_2,Rxy));
&colore{darkorange}{i12 :} &colore{darkorchid}{time} &colore{bleudefrance}{sparseResultant}(f_0,f_1,f_2) == 
           -&colore{bleudefrance}{sparseDiscriminant}(f'_0 + y_1*f'_1 + y_2*f'_2)
&colore{commentoutput}{      -- used 0.746274 seconds}
&colore{circOut}{o12 =} &colore{output}{true}
\end{Verbatim}
} 
\end{tcolorbox}
\end{example}

A derived function of \texttt{sparseDiscriminant} is \texttt{determinant} (or simply \texttt{det}), which calculates determinants of multidimensional matrices (see subsection~\ref{sec dets}).
However for this last one, more specialized algorithms are also available and automatically applied.
\begin{example}
Here we calculate the determinant of a generic 
four-dimensional
matrix 
of shape $2\times 2\times 2\times 2$ (see also \cite{2x2x2x2}). This polynomial is homogeneous of degree $24$ in 
the $16$ variable entries of the matrix and it has $2,894,276$ terms. 
The approach for this calculation is to apply Theorem~\ref{sch} recursively.
The calculation time  
is about $10$ minutes,
but it takes much less time 
if we specialize the entries of the matrix to be random numbers.
\begin{tcolorbox}[breakable=true,boxrule=0pt,opacityback=0.0,enhanced jigsaw]
{\footnotesize
\begin{Verbatim}[commandchars=&!$]
&colore!darkorange$!i13 :$ M = &colore!bleudefrance$!genericMultidimensionalMatrix$ {2,2,2,2}
&colore!circOut$!o13 =$ &colore!output$!{{{{a       , a       }, {a       , a       }}, {{a       , a       }, ...$
&colore!output$!           0,0,0,0   0,0,0,1     0,0,1,0   0,0,1,1       0,1,0,0   0,1,0,1   ...$
&colore!circOut$!o13 :$ &colore!output$!4-dimensional matrix of shape 2 x 2 x 2 x 2 over ZZ[a       , a       , ...$
&colore!output$!                                                           0,0,0,0   0,0,0,1  ...$
&colore!darkorange$!i14 :$ &colore!darkorchid$!time$ D = &colore!bleudefrance$!det$ M;
&colore!commentoutput$!      -- used 634.773 seconds$
&colore!darkorange$!i15 :$ (&colore!airforceblue$!first degree$ D, # &colore!airforceblue$!terms$ D)
&colore!circOut$!o15 =$ &colore!output$!(24, 2894276)$
\end{Verbatim}
} 
\end{tcolorbox}
\end{example}

\begin{example} 
Here we take $A$ and $B$ to be random matrices of shapes $2\times2\times 2\times4$ and $4\times 2\times 5$, respectively.
We calculate the convolution $A\ast B$ (see \cite[p.~449]{gelfandbook}), which is a matrix of shape $2\times 2\times 2 \times 2 \times 5$.
Then we verify a formula 
proved in \cite{DIONISI200387} for $\det(A\ast B)$,
which generalizes the Cauchy-Binet formula in the multidimensional case.
 The approach for the calculation of the determinant of
 shape $4\times2\times 5$ is using Proposition~\ref{boundary det},
 while the determinants of shapes $2\times2\times2\times 4$ and $2\times2\times2\times2\times 5$ 
 are calculated using Remark~\ref{remBoundaryDet}.
 The calculation time is less than one second.
\begin{tcolorbox}[breakable=true,boxrule=0pt,opacityback=0.0,enhanced jigsaw]
{\footnotesize
\begin{Verbatim}[commandchars=&!$]
&colore!darkorange$!i16 :$ K = ZZ/33331;
&colore!darkorange$!i17 :$ A = &colore!bleudefrance$!randomMultidimensionalMatrix$({2,2,2,4},&colore!airforceblue$!CoefficientRing$=>K);
&colore!circOut$!o17 :$ &colore!output$!4-dimensional matrix of shape 2 x 2 x 2 x 4 over K$
&colore!darkorange$!i18 :$ B = &colore!bleudefrance$!randomMultidimensionalMatrix$({4,2,5},&colore!airforceblue$!CoefficientRing$=>K);
&colore!circOut$!o18 :$ &colore!output$!3-dimensional matrix of shape 4 x 2 x 5 over K$
&colore!darkorange$!i19 :$ &colore!darkorchid$!time$ &colore!bleudefrance$!det$(A * B) == (&colore!bleudefrance$!det$ A)^5 * (&colore!bleudefrance$!det$ B)^6
&colore!commentoutput$!      -- used 0.535271 seconds$
&colore!circOut$!o19 =$ &colore!output$!true$
\end{Verbatim}
} 
\end{tcolorbox}
\end{example}


\begin{thebibliography}{JMSW09}

\bibitem[Cay45]{cayley_2009}
A.~Cayley, \emph{On the theory of linear transformations}, Cambridge Math. J.
  \textbf{4} (1845), 1--16.

\bibitem[CDS98]{CattaniDickenstein1998}
E.~Cattani, A.~Dickenstein, and B.~Sturmfels, \emph{Residues and resultants},
  J. Math. Sci. Univ. Tokyo \textbf{5} (1998), 119--148.

\bibitem[CE00]{CaEmi2000}
J.~F. Canny and I.~Z. Emiris, \emph{A subdivision-based algorithm for the
  sparse resultant}, J. {ACM} \textbf{47} (2000), 417--451.

\bibitem[CLO05]{coxbook}
D.~Cox, J.~Little, and D.~O'Shea, \emph{Using algebraic geometry}, second ed.,
  Grad. Texts in Math., vol. 185, Springer, 2005.

\bibitem[D'A02]{DAndrea2001MacaulaySF}
C.~D'Andrea, \emph{Macaulay style formulas for sparse resultants}, Trans. Amer.
  Math. Soc. \textbf{354} (2002), no.~7, 2595--2629.

\bibitem[Dix08]{Dixon1908}
A.~L. Dixon, \emph{The eliminant of three quantics in two independent
  variables}, Proc. Lond. Math. Soc. \textbf{7} (1908), 49--69.

\bibitem[DO03]{DIONISI200387}
C.~Dionisi and G.~Ottaviani, \emph{The {B}inet-{C}auchy theorem for the
  hyperdeterminant of boundary format multi-dimensional matrices}, J. Algebra
  \textbf{259} (2003), no.~1, 87--94.

\bibitem[DS15]{PoissonSombraAndrea}
C.~D'Andrea and M.~Sombra, \emph{A {P}oisson formula for the sparse resultant},
  Proc. Lond. Math. Soc. \textbf{110} (2015), no.~4, 932--964.

\bibitem[EM99]{EMIRIS19993}
I.~Z. Emiris and B.~Mourrain, \emph{Matrices in elimination theory}, J.
  Symbolic Comput. \textbf{28} (1999), no.~1, 3--44.

\bibitem[GKZ92]{GELFAND1992226}
I.~M. Gelfand, M.~M. Kapranov, and A.~V. Zelevinsky, \emph{Hyperdeterminants},
  Adv. Math. \textbf{96} (1992), no.~2, 226--263.

\bibitem[GKZ94]{gelfandbook}
\bysame, \emph{Discriminants, resultants, and multidimensional determinants},
  Mathematics: Theory \& Applications, Birkh{\"a}user Boston, 1994, reprinted
  in 2008.

\bibitem[GS21]{macaulay2}
D.~R. Grayson and M.~E. Stillman, \emph{{\sc Macaulay2} --- {A} software system
  for research in algebraic geometry (version 1.18)}, home page:
  \url{http://www.math.uiuc.edu/Macaulay2/}, 2021.

\bibitem[HSYY08]{2x2x2x2}
P.~Huggins, B.~Sturmfels, J.~Yu, and D.~S. Yuster, \emph{The hyperdeterminant
  and triangulations of the $4$-cube}, Math. Comp. \textbf{77} (2008), no.~263,
  1653--1679.

\bibitem[JKSS04]{JeronimoSabia2004}
G.~Jeronimo, T.~Krick, J.~Sabia, and M.~Sombra, \emph{The computational
  complexity of the {C}how form}, Found. Comput. Math. \textbf{4} (2004),
  41--117.

\bibitem[JMSW09]{Matera2009Jeronimo}
G.~Jeronimo, G.~Matera, P.~Solern\'o, and A.~Waissbein, \emph{Deformation
  techniques for sparse systems}, Found. Comput. Math. \textbf{9} (2009),
  1--50.

\bibitem[Mac02]{macaulaypaper1902}
F.~Macaulay, \emph{On some formulas in elimination}, Proc. Lond. Math. Soc.
  \textbf{3} (1902), 3--27.

\bibitem[Ott13]{OttSurvey}
G.~Ottaviani, \emph{Introduction to the hyperdeterminant and to the rank of
  multidimensional matrices}, Commutative Algebra: Expository Papers Dedicated
  to {D.} {E}isenbud on the Occasion of His 65th Birthday (2013).

\bibitem[Sch52]{bhlpart219634}
L.~Schl{\"a}fli, \emph{{\"U}ber die resultante eines systemes mehrerer
  algebraischer gleichungen. {E}in beitrag zur theorie der elimination},
  Denkschr der Kaiserlichen Akad. der Wiss, Math-Naturwiss. Classe \textbf{4}
  (1852), 1--74.

\bibitem[Sta18]{packageResultants}
G.~Staglian\`o, \emph{A package for computations with classical resultants}, J.
  Softw. Alg. Geom. \textbf{8} (2018), no.~1, 21--30.

\bibitem[Stu93]{SturmfelsElimination}
B.~Sturmfels, \emph{Sparse elimination theory}, Computational algebraic
  geometry and commutative algebra (D.~Eisenbud and L.~Robbiano, eds.),
  Cambridge Univ. Press, Cambrige, 1993, pp.~264--298.

\bibitem[Stu94]{NewtonSturm}
\bysame, \emph{On the {N}ewton polytope of the resultant}, J. Algebraic Combin.
  \textbf{3} (1994), 207--236.

\bibitem[Stu02]{Sturmfels02solvingsystems}
\bysame, \emph{Solving systems of polynomial equations}, {CBMS} Regional Conf.
  Ser. in Math., vol.~97, Amer. Math. Soc., 2002.

\bibitem[WZ96]{AIF19964635910}
J.~Weyman and A.~Zelevinsky, \emph{Singularities of hyperdeterminants}, Ann.
  Inst. Fourier (Grenoble) \textbf{46} (1996), no.~3, 591--644.

\end{thebibliography}

\providecommand{\bysame}{\leavevmode\hbox to3em{\hrulefill}\thinspace}
\providecommand{\MR}{\relax\ifhmode\unskip\space\fi MR }
\providecommand{\MRhref}[2]{%
  \href{http://www.ams.org/mathscinet-getitem?mr=#1}{#2}
}
\providecommand{\href}[2]{#2}

\end{document}